\documentclass{article}
\usepackage{amsmath} 
\usepackage{amsthm}  
\usepackage{cite}

\def\firstpage{1}

\title{\Large \bfseries A Laplace transform approach to $C$-semigroups on a  $\mathcal{T}_{\varepsilon, \lambda}$-complete random normed module}
\author{\large Xia Zhang$^1$\qquad Leilei Wei$^{2}$\qquad Ming Liu$^3$}
\date{}

\begin{document}

\maketitle

\thispagestyle{empty}
\renewcommand{\thefootnote}{\fnsymbol{footnote}}

\footnotetext{\hspace*{-5mm} \begin{tabular}{@{}r@{}p{13.4cm}@{}}
& Manuscript received  \\
$^1$ & School of Mathematical Sciences, TianGong University, Tianjin 300387, P.R.China\\
&{E-mail: zhangxia@tiangong.edu.cn} \\
$^{2}$ & School of Mathematical Sciences, TianGong University, Tianjin 300387, P.R.China\\
&{E-mail: 2220110138@tiangong.edu.cn} \\
$^{3}$ & Corresponding author. School of Mathematical Sciences, TianGong University, Tia-\\
& njin 300387, P.R.China\\
&{E-mail: liuming@tiangong.edu.cn} \\
 $^{\ast}$ & Project supported by
  the National Natural Science Foundation of China (Grant No.\\& 12171361)
   and the Humanities and Social Science Foundation of Ministry of Educati-\\
   &on (Grant No. 20YJC790174).
\end{tabular}}

\renewcommand{\thefootnote}{\arabic{footnote}}

\begin{abstract}
In this paper, we first introduce the notion of the Laplace transform for an abstract-valued function from $[0, \infty)$ to a $\mathcal{T}_{\varepsilon, \lambda}$-complete  random normed module $S$. Then, combining respective advantages of the  $(\varepsilon, \lambda)$-topology and the locally $L^0$-convex topology on $S$,
we prove the differentiability, Post-Widder inversion formula and uniqueness of such a Laplace transform. Second, based on the above work, we
establish the Hille-Yosida type theorem for an exponentially bounded $C$-semigroup on  $(S, \mathcal{T}_{\varepsilon, \lambda})$ for the case that $R(C)$
is  nondense in  $(S, \mathcal{T}_{\varepsilon, \lambda})$, which extends and improves several important results. Finally, we also apply such a Laplace transform to abstract Cauchy problems in the random setting.

\vskip 4.5mm

\noindent
\begin{tabular}{@{}l@{ }p{10.1cm}}
\textbf{Keywords} &
Random normed modules, Laplace transforms, $C$-semigroups, \\
&Abstract Cauchy problems
\end{tabular}

\noindent
\textbf{2000 MR Subject Classification}
46H25,45R05

\end{abstract}

\renewcommand{\baselinestretch}{1.05}\selectfont
\setlength{\parindent}{1.5em}
\setcounter{section}{0}

\section{Introduction} \label{section1}

\quad\;The approach of Laplace transforms is closely related to the theory of operator semigroups and abstract Cauchy problems \cite{Arendt, deLaubenfels1}. In particular, the establishment of the operator semigroup theory also marks the beginning of the study of Laplace transforms for abstract-valued functions from \([0, \infty)\) to Banach spaces.
The key figures of this stage included mathematicians such as Hille, Yosida and Phillips, who established the Laplace transform of such abstract-valued functions and linked it to the theory of operator semigroups \cite{Hille48, Phillips1954, Yosida1948}. For example, the Hille-Yosida theorem, as the cornerstone of operator semigroup theory, can be regarded as the embodiment of such abstract-valued Laplace transforms on the resolvent of a specific operator (the infinitesimal generator). Therefore, deepening the research on Laplace transform theory is an important way to promote the development of operator semigroup theory.

The notion of a random normed module is a random generalization of that of an ordinary normed space and is also one of the central frameworks of random functional analysis. In the past more than thirty years, thanks to Guo's pioneering contributions \cite{Guo01,Guo10,Guo20241}, the theory of random normed modules has obtained a systematic and deep development \cite{TT,Guo2024,Guo0912,Guozhao,Guozhao15,Thang19,Wu2013,WuGuo2015,WuGuoLong2022,WuZengZhao2022,Zhangliuguo2}, and has already been successfully applied to several important related  fields such as random operator equations \cite{GuoZhangwuyang17,LIU2,LIU1}, dynamic risk measures \cite{GuoZhaozeng14} and nonsmooth differential geometry on metric measure spaces \cite{Guo2024,E2019,E}.
It is known that the norm of a normed space can induce only one strong topology, i.e., the norm topology, whereas the $L^{0}$-norm of a random normed module can induce two kinds of strong topologies~--~the $(\varepsilon, \lambda)$-topology (denoted by $\mathcal{T}_{\varepsilon, \lambda}$) and the locally $L^{0}$-convex topology (denoted by $\mathcal{T}_c$), which are both frequently employed in the development of random normed modules. Moreover, such two kinds of topologies possess their apparent advantages and disadvantages, respectively. For example, the $(\varepsilon, \lambda)$-topology is very natural from the viewpoint of probability theory, which makes  Guo et al. successfully establish a hyperplane separation theorem between a point and a $\mathcal{T}_{\varepsilon, \lambda}$-closed $L^{0}$-convex subset in \cite{Guo0912}, a $(\varepsilon, \lambda)$-topological version of Fenchel-Moreau duality theorem for a lower semicontinuous $L^{0}$-valued function on random locally convex modules in \cite{Guozhao15}, and the fundamental theorem of calculus in $\mathcal{T}_{\varepsilon, \lambda}$-complete random normed modules \cite{Guozhang12}. However, the $(\varepsilon, \lambda)$-topology is too weak to ensure that a random normed module admits a proper  $\mathcal{T}_{\varepsilon, \lambda}$-open $L^{0}$-convex subset, whereas the locally $L^{0}$-convex topology can make a random normed module admit proper $\mathcal{T}_{c}$-open $L^{0}$-convex  subsets  so that many propositions in which some $L^0$-convex subsets are required to have nonempty interiors can be established under such a topology in \cite{GuoZhaozeng14}, but it is too strong to ensure that many classical propositions to hold in a natural way. A major breakthrough during this period was that Guo established an intrinsic connection between the two kinds of topologies by introducing the concept of $\sigma$-stability (the countable concatenation property) \cite{Guo10}.

It is well known that the Laplace transform for abstract-valued functions from \([0, \infty)\) to Banach spaces has been widely applied to differential equations and operator semigroups. However, this transform universally fails to serve for the theory of random normed modules, which motivates us to introduce the notion of Laplace transforms in the random setting. In this paper, we just simultaneously use the advantages of the two kinds of topologies to complete the proofs of the differentiability, Post-Widder inversion formula, and uniqueness of the Laplace transform in complete random normed modules. The difficult point in this paper lies in the proof of Theorem 3.1.
It is known from \cite{Guozhao} that a $\mathcal{T}_c$-complete random normed module with base $(\Omega, \mathcal{F}, P)$ is a totally disconnected topological space when $(\Omega, \mathcal{F}, P)$ is nonatomic, at which time a continuous function from a finite real interval to a $\mathcal{T}_c$-complete random normed module can only be a constant function. This fact also shows that it makes no sense to define the Riemann integral for such continuous functions. So Guo and Zhang turned to define the Riemann integral for abstract-valued functions from a finite real closed interval  to a $\mathcal{T}_{\varepsilon, \lambda}$-complete random normed module in \cite{Guozhang12}, and further established the fundamental theorem of calculus in such random setting (see Proposition 2.1 in this paper), which played an important role in the subsequent theory of operator semigroups in $\mathcal{T}_{\varepsilon, \lambda}$-complete random normed modules.
Therefore, to define the Laplace transform in the random setting, we have to use the $(\varepsilon, \lambda)$-topology. On the other hand, since the $(\varepsilon, \lambda)$-topology is essentially not locally convex in general, there may be no $\mathcal{T}_{\varepsilon, \lambda}$-open $L^0$-convex sets in random normed  modules. Thus, to establish the differentiability of such Laplace transforms, we are forced to turn to $\mathcal{T}_c$-open subsets, which necessitates the simultaneous use of the locally $L^0$-convex topology in $\mathcal{T}_{\varepsilon, \lambda}$-complete random normed modules.

The theory of $C$-semigroups on a Banach space $X$ is an important generalization of $C_{0}$-semigroups on $X$ \cite{Davies,Miyadera86} and provides power tools for dealing with ill-posed abstract Cauchy problems  \cite{deLaubenfels1}.
Recently, the study of operator semigroups in random normed modules has also obtained some progress \cite{Thang19,Zhangliu13,Zhangliuguo2}. In particular, in 2024, Son, Thang and Oanh first studied the exponentially bounded $C$-semigroups and the Cauchy initial value problems in $\mathcal{T}_{\varepsilon, \lambda}$-complete random normed modules \cite{Son}.  Subsequently, Zhang and Liu studied the abstract Cauchy problem for evolution equations in $\mathcal{T}_{\varepsilon, \lambda}$-complete random normed modules \cite{Zhangliu25}.
Recently, Zhang, Wei and Liu established some relations among $C$-existence families, $C$-semigroups and their associated abstract Cauchy problems in $\mathcal{T}_{\varepsilon, \lambda}$-complete random normed modules \cite{Zhangliuwei}. 
The aim of this paper is to continue to study the theory of exponentially bounded $C$-semigroups in the random setting, specially, we will first introduce the notion of the Laplace transform for an abstract-valued function from $[0, \infty)$ to a $\mathcal{T}_{\varepsilon, \lambda}$-complete random normed module and further apply such Laplace transforms to exponentially bounded  $C$-semigroups on a  $\mathcal{T}_{\varepsilon, \lambda}$-complete random normed module.

This paper contains five sections: in Section 2, we will present some preliminaries.  In Section 3, we will introduce the notion of the Laplace transform for an abstract-valued function from $[0, \infty)$ to a $\mathcal{T}_{\varepsilon, \lambda}$-complete random normed module $S$ and present some basic results peculiar to such Laplace transforms; then, based on the above work, in Section 4, we will establish
the Hille-Yosida type theorem for an exponentially bounded $C$-semigroup on $(S, \mathcal{T}_{\varepsilon, \lambda})$ for the case that $R(C)$
is  nondense in  $(S, \mathcal{T}_{\varepsilon, \lambda})$. In Section 5,  we will apply such Laplace transforms to abstract Cauchy problems in the random setting.

\section{Preliminaries}

\quad\;In this paper, we start with some notations. N denotes the set of positive integers, K the real scalar field R or the complex scalar field C, and $(\Omega, \mathcal{F}, P)$ a given probability space. Moreover, $L^0(\mathcal{F}, K)$ denotes the algebra of equivalence classes of K-valued $\mathcal{F}$-measurable random variables on $\Omega$, and $\bar{L}^0(\mathcal{F}, R)$ the set of equivalence classes of extended real-valued $\mathcal{F}$-measurable random variables on $\Omega$. Clearly, it is known from \cite{D} that $\bar{L}^0(\mathcal{F}, R)$ is a complete lattice under the partial order $\leq:$ $f \leq g$ if and only if $f^0(\omega) \leq g^0(\omega)$ for almost all $\omega$ in $\Omega$, where $f^0$ and $g^0$ are arbitrarily chosen representatives of $f$ and $g$ in $\bar{L}^0(\mathcal{F}, R)$.
Besides, for any $G \subset \bar{L}^0(\mathcal{F}, R)$ and $G \neq \emptyset$, let $\bigvee G$ and $\bigwedge G$ denote the supremum and infimum of $G$, respectively, then there are two sequences $\left\{\xi_n, n \in N\right\}$ and $\left\{\eta_n, n \in N\right\}$ in $G$ satisfying that $\bigvee_{n \geq 1} \xi_n=\bigvee G$ and $\bigwedge_{n \geq 1}$ $\eta_n=\bigwedge G$. Further, $L^0(\mathcal{F}, R)$, as a sublattice of $\bar{L}^0(\mathcal{F}, R)$, is complete in the sense that each subset with an upper bound has a supremum.

 As usual, $I_D$ denotes the characteristic function of $D$ for any $D \in \mathcal{F}$ and $\tilde{I}_D$ denotes the equivalence class of $I_D$. For any $f, g \in \bar{L}^0(\mathcal{F}, R)$, $f>g$ means  $f \geq g$ and $f \neq g$, and for any $A \in \mathcal{F}$, $f>g$ on $A$ means $f^0(\omega)>g^0(\omega)$ for almost all $\omega \in A$, where $f^0$ and $g^0$ are arbitrarily chosen representatives of $f$ and $g$, respectively.  Besides, let $D=\left\{\omega \in \Omega \mid f^0(\omega)>\right.$ $\left.g^0(\omega)\right\}$, then we always use $[f>g]$ for the equivalence class of $D$ and often write $I_{[f>g]}$ for $\tilde{I}_D$, one can also understand such notations as $I_{[f \leq g]}, I_{[f \neq g]}$ and $I_{[f=g]}$.

 In particular, we denote  $L_{+}^0(\mathcal{F})=\{\xi \in L^0(\mathcal{F}, R)$ $ \mid \xi \geq 0\}$ and $L_{++}^0(\mathcal{F})=\{\xi \in$ $L^0(\mathcal{F}, R) \mid \xi>0$ on  $\Omega\}$.

\textbf{Definition 2.1} \textrm{(see \cite{Guo01})}
An ordered pair $(S,\|\cdot\|)$ is called a random normed module (briefly, an $RN$ module) over K with base $(\Omega, \mathcal{F}, P)$ if $S$ is  a left module over the algebra $L^0(\mathcal{F}, K)$ and $\|\cdot\|$ is a mapping from $S$ to $L_{+}^0(\mathcal{F})$ such that the following three axioms are satisfied:

(1) $\|x\|=0$ $\Rightarrow$ $x=\theta$ (the null in $S$);

(2) $\|\zeta x\|=|\zeta| \cdot\|x\|,$ $\forall\zeta \in L^0(\mathcal{F}, K)$ and $x \in S$;

(3) $\|x+y\| \leq\|x\|+\|y\|,$ $\forall x, y\in S$.

As usual, $\|\cdot\|$ is called an $L^0$-norm on $S$ and $\|x\|$ is called the $L^0$-norm of $x$ in $S$.

For example,  $\left(L^0(\mathcal{F}, K),|\cdot|\right)$ is a nontrivial $R N$ module, where $|\cdot|$ denotes the absolute value function.

Let $(S_1,\|\cdot\|_1)$ and $(S_2,\|\cdot\|_2)$ be two $R N$ modules over K with base $(\Omega, \mathcal{F}, P)$. A linear operator $T$ from $S_1$ to $S_2$ is said to be almost surely bounded if there is an $\eta \in L_{+}^0(\mathcal{F})$ satisfying $\|T z\|_2 \leq \eta \cdot\|z\|_1$ for any $z \in S_1$. Denote by $B\left(S_1, S_2\right)$ the linear space of almost surely bounded linear operators from $S_1$ to $S_2$, and define a mapping $\|\cdot\|: B\left(S_1, S_2\right) \rightarrow L_{+}^0(\mathcal{F})$ by $\|T\|:=\bigwedge\left\{\eta \in L_{+}^0(\mathcal{F}) \mid\|T z\|_2 \leq \eta \cdot\|z\|_1\right.$ for any $\left.z \in S_1\right\}$ for any $T \in B\left(S_1, S_2\right)$, then one can obtain that $\left(B\left(S_1, S_2\right),\|\cdot\|\right)$ is still an $R N$ module.
It is known from \cite{TT} that $T$ is almost surely bounded if and only if $T$ is a continuous module homomorphism from $S_1$ to $S_2$. Besides, if
$S_1=S_2=S$, then $B(S_1, S_2)$ is abbreviated as $B(S)$.

\textbf{Definition 2.2} \textrm{(see \cite{GuoZhaozeng14})}
Let $(S,\|\cdot\|)$ be an $R N$ module over K with base $(\Omega, \mathcal{F}, P)$. For any $\varepsilon>0$ and $0<\lambda<1$, let $U_\theta(\varepsilon, \lambda)=\{z \in S \mid P\{\omega \in \Omega \mid\|z\|(\omega)<\varepsilon\}>1-\lambda\}$ and $\mathcal{U}_\theta=\left\{U_\theta(\varepsilon, \lambda) \mid \varepsilon>0,0<\lambda<1\right\}$, then $\mathcal{U}_\theta$ forms a local base at $\theta$ of some metrizable linear topology, called the $(\varepsilon, \lambda)$-topology, denoted by $\mathcal{T}_{\varepsilon, \lambda}$.

It is clear that a sequence $\{z_{n},n\in N\}\subset S$ converges to $z_{0}\in S$ in the $(\varepsilon,\lambda)$-topology if and only if the sequence $\{\|z_{n}-z_{0}\|,n\in N\}\subset L^0(\mathcal{F}, R)$ converges to 0 in probability $P$.  In particular, the
$(\varepsilon,\lambda)$-topology for $(L^0(\mathcal{F}, K),|\cdot|)$ is exactly the one of convergence in  probability $P$.

 \textbf{Definition 2.3} \textrm{(see \cite{Guo10})}
Let $(S,\|\cdot\|)$ be an $R N$ module over K with base $(\Omega, \mathcal{F}, P)$. For any $ \varepsilon \in L_{++}^0(\mathcal{F})$, let $U_\theta(\varepsilon)=\{x \in S~|~\|x\|<\varepsilon$ on $\Omega\}$ and $\mathcal{U}_\theta=\left\{U_\theta(\varepsilon)~|~ \varepsilon \in L_{++}^0(\mathcal{F})\right\}$.  A set $G \subset S$ is called $\mathcal{T}_c$-open if for every $x \in G$ there exists some $U_\theta (\varepsilon)\in \mathcal{U}_\theta$ such that $x+U_\theta (\varepsilon)\subset G$. Let $\mathcal{T}_c$ be the family of $\mathcal{T}_c$-open subsets, then $\mathcal{T}_c$ is a Hausdorff topology on $S$, called the locally $L^0$-convex topology.

It is worthnoting that the locally $L^{0}$-convex topology is not necessarily linear and such a topology  for the algebra $L^0(\mathcal{F}, K)$ is only a topological ring in general. Moreover, $\left(S, \mathcal{T}_c\right)$ is a topological module over the topological ring $\left(L^0(\mathcal{F}, K), \mathcal{T}_c\right)$ and $\left\{U_\theta(\varepsilon)~|~ \varepsilon \in L_{++}^0(\mathcal{F})\right\}$ is just a local base at $\theta$ of $\mathcal{T}_c$.
In general, $\mathcal{T}_c$ is not metrizable but it can ensure that some $L^0$-convex subsets have nonempty interiors.
A net $\left\{x_\alpha, \alpha \in \wedge\right\}$ in $S$ converges in the locally $L^0$-convex topology to $x \in S$ if and only if $\left\{\left\|x_\alpha-x\right\|\right.$, $\alpha \in \wedge\}$ converges in the locally $L^0$-convex topology of $L^0(\mathcal{F}, K)$ to $\theta$.

Let $[s, t]$ be a finite closed real interval and $(S,\|\cdot\|)$ an $R N$ module over K with base $(\Omega, \mathcal{F}, P)$. A function
$g:[s, t] \rightarrow S$ is said to be $L^0$-Lipschitz on $[s, t]$ if there is an $\eta \in L_{+}^0(\mathcal{F})$ satisfying that $\left\|g\left(s_1\right)-g\left(s_2\right)\right\| \leq \eta\left|s_1-s_2\right|$ for any $s_1, s_2 \in[s, t]$. Furthermore, a function $g:[0, \infty) \rightarrow S$ is said to be locally $L^0$-Lipschitz if for any $L>0$, there is a $\zeta_L \in L_{+}^0(\mathcal{F})$ satisfying $\left\|g\left(s_1\right)-g\left(s_2\right)\right\| \leq \zeta_L\left|s_1-s_2\right|$ for any $s_1, s_2 \in[0, L]$.

\textbf{Proposition 2.1} \textrm{(see \cite{Guozhang12})}
Suppose that $S$ is a $\mathcal{T}_{\varepsilon, \lambda}$-complete $R N$ module and a function $g:[s, t] \rightarrow (S, \mathcal{T}_{\varepsilon, \lambda})$ is continuously differentiable. If $g$ is $L^0$-Lipschitz on $[s, t]$, then $g^{\prime}$ is Riemann integrable on $[s, t]$ and
$$
\int_s^t g^{\prime}(u) du=g(t)-g(s).
$$

As we all know, a continuous function from $[s, t]$ to a Banach space is automatically bounded, but a continuous function from $[s, t]$ to a $\mathcal{T}_{\varepsilon, \lambda}$-complete $R N$ module $S$ may not be almost surely  bounded. Fortunately, a sufficient condition for a continuous function  to be Riemann integrable has been given, that is, if $g:[s, t] \rightarrow (S, \mathcal{T}_{\varepsilon, \lambda})$ is a continuous function satisfying that $\bigvee_{u \in[s, t]}\|g(u)\|$ belongs to $L_{+}^0(\mathcal{F})$, then $g$ is Riemann integrable. Based on this fact, Propositions 2.2 and 2.3 below hold.

\textbf{Proposition 2.2} \textrm{(see \cite{Guozhang12})}
Suppose that $S$ is a $\mathcal{T}_{\varepsilon, \lambda}$-complete $R N$ module and a function $g:[s, t] \rightarrow (S, \mathcal{T}_{\varepsilon, \lambda})$ is continuous. If $\bigvee_{u \in[s, t]}\|g(u)\| \in L_{+}^0(\mathcal{F})$, then the following statements hold.

(1) $\|\int_s^t g(u) d u\| \leq \int_s^t\|g(u)\| d u$;

(2) Let $G(l)=\int_s^l g(u) d u$ for any $l \in[s, t]$, then $G$ is $\mathcal{T}_{\varepsilon, \lambda}$-differentiable on $[s, t]$ and $G^{\prime}(l)=g(l)$.

\textbf{Proposition 2.3} \textrm{(see \cite{Zhangliu13})}
Suppose that $g:[s, t] \rightarrow (L^0(\mathcal{F}, R), \mathcal{T}_{\varepsilon, \lambda})$ is a continuous function  satisfying $\bigvee_{u \in[s, t]}|g(u)| \in L^1(\mathcal{F}, R)$, where $L^1(\mathcal{F}, R)=\{f \in L^0(\mathcal{F}, R)$ $\mid \int_{\Omega}|f| d P<\infty\}$. Then
$$
\int_{\Omega}\left[\int_s^t g(u) d u\right] d P=\int_s^t\left[\int_{\Omega} g(u) d P\right] d u .
$$

\section{ The Laplace transform for an abstract-valued function from $[0, \infty)$ to a $\mathcal{T}_{\varepsilon, \lambda}$-complete  $R N$ module}\label{section3}
\quad\;This section is devoted to introducing   the notion of  the Laplace transform for an abstract-valued function from $[0, \infty)$ to  a $\mathcal{T}_{\varepsilon, \lambda}$-complete $R N$ module, and further establishing  its differentiability, the Post-Widder inversion formula and uniqueness theorems, which play a crucial role in the proof
of Theorem 4.1.

For any $\xi \in L^0(\mathcal{F}, R)$, let $G_\xi=\left\{\eta \in L^0(\mathcal{F}, R)\mid \eta> \xi\; on\;  \Omega\right\}$, then $G_\xi$ is a $\mathcal{T}_c$-open subset of $L^0(\mathcal{F}, R)$. Clearly, $G_0=\left\{\eta \in L^0(\mathcal{F}, R)\mid \eta> 0\; on\;  \Omega\right\}$.

\textbf{Definition 3.1}
Let $S$ be a $\mathcal{T}_{\varepsilon, \lambda}$-complete $R N$ module and  $h:[0, \infty) \rightarrow (S, \mathcal{T}_{\varepsilon, \lambda})$ a continuous function satisfying $\|h(s)\| \leq M e^{\xi s}$ for any $s\geq 0$ and some $M \in L_{+}^0(\mathcal{F})$, $\xi \in L^0(\mathcal{F}, R)$, then the Laplace transform of $h$ is given by
$$
H(\eta)=\int_0^\infty e^{-\eta s}h(s)ds
$$
for any $\eta \in G_\xi$.

In fact, for any $\eta \in G_\xi$,
\allowdisplaybreaks[4] 
\begin{align*}
\|H(\eta)\| & =\|\int_0^{\infty} e^{-\eta s}h(s)d s\| \\
& \leq \int_0^{\infty}e^{-\eta s}\|h(s)\| ds \\
& \leq \frac{M}{\eta-\xi},
\end{align*}
thus the Laplace transform is well-defined.

\textbf{Definition 3.2}
Let $S$ be an $R N$ module, $G$ a $\mathcal{T}_c$-open subset of $L^0(\mathcal{F}, R)$ and $H:G\rightarrow S$ a  function.
For any $x_0\in G$, if there exists an  $I \in S$  such that  $\frac{H\left(x_0+h\right)-H\left(x_0\right)}{h}-I$  converges  to  $\theta$ in $\mathcal{T}_c$
when  $h \in L^0(\mathcal{F}, R)$ such that  $|h| \in L_{++}^0(\mathcal{F})$ and  $h$ converges to  $0$ in $\mathcal{T}_c$, then $H$ is said to be $\mathcal{T}_c$-differentiable at $x_0$ and $I$ is called the $\mathcal{T}_c$-derivative of $H$ at $x_0$, denoted by $\mathcal{T}_c\text{-}H^{\prime}(x_0)$ or $\mathcal{T}_c\text{-}\frac{d H(x_0)}{dx_0}$.

\textbf{Lemma 3.1}
For any $s\geq 0$ and $h \in  L^0(\mathcal{F}, R)$ such that $|h| \in L_{++}^0(\mathcal{F})$, there exists a unique $t\in L_{++}^0(\mathcal{F})$ with $t<s$ on $\Omega$ such that $$e^{-h s}=1-h s+\frac{(h s)^2}{2} e^{-h t}.$$

\textbf{Proof}
Let $h^0$ be an arbitarily chosen representative of $h$ such that $0<h^{0}(\omega)<\infty$ for any $\omega \in \Omega$. Now, for any $\omega \in \Omega$, by applying Maclaurin's formula with the lagrange remainder to the function $f_\omega:[0, s] \rightarrow [0, \infty)$ defined by $f_\omega(r)=e^{-h^0(\omega) r}$, there exists some $t^0(\omega) \in(0, s)$ such that $f_\omega(s)=1-h^0(\omega) s+\frac{\left(h^0(\omega) s\right)^2}{2} e^{-h^0(\omega) t^0(\omega)}$. In fact, $t^0(\omega)$ is unique for any $\omega \in \Omega$, and $t^0(\omega)=\frac{-\ln \left(\left|\frac{2(e^{-h^0(\omega) s}-1+h^0(\omega) s) }{\left(h^0(\omega) s\right)^2}\right|\right)}{ h^0(\omega)}$. Clearly, $t^0: \Omega \rightarrow(0, s)$ is $\mathcal{F}$-measurable. Let $t$ be the equivalence class of $t^0$, then $t$ is desired.
\endproof

\textbf{Lemma 3.2} (Differentiability)
Suppose that $S$ is a $\mathcal{T}_{\varepsilon, \lambda}$-complete $R N$ module and $h:[0, \infty) \rightarrow (S, \mathcal{T}_{\varepsilon, \lambda})$ is a continuous function satisfying  $\|h(s)\| \leq M e^{\xi s}$ for any $s\geq 0$ and some $M \in L_{+}^0(\mathcal{F})$, $\xi \in L^0(\mathcal{F}, R)$. For any $\eta \in G_\xi$,
let $H(\eta)$ be the Laplace transform of $h$, then
$$
\mathcal{T}_c\text{-}H^{(k)}(\eta)=\int_0^{\infty} e^{-\eta s}(-s)^k h(s)d s
$$
for any $k \in N$.

\textbf{Proof}
For any $\eta \in G_\xi$, it is clear that there exists a $\gamma \in L_{++}^0(\mathcal{F})$ such that $\eta> \xi+\gamma$ on  $\Omega$. Let $\bar{G}=\left\{\zeta \in L^0(\mathcal{F}, R)\mid |\zeta-\eta| < \eta-\xi-\gamma\; on\;  \Omega\right\}$, then $\bar{G}$ is a $\mathcal{T}_c$-open subset of $G_\xi$ and $\eta \in \bar{G}$.  Further, set $$\bar{\bar{G}}=\left\{h \in L^0(\mathcal{F}, R)\mid |h| \in L_{++}^0(\mathcal{F})\; and\; \eta+h \in \bar{G}\right\},$$
then it is easy to check that  $|h| < \eta- \xi-\gamma$ on  $\Omega$ for any $h \in \bar{\bar{G}}$.

Without loss of generality, we can assume that $k=1$.
For any $s\geq 0$ and $h \in \bar{\bar{G}}$, according to Lemma 3.1, there exists a unique $t\in L_{++}^0(\mathcal{F})$ with $t<s$ on $\Omega$ such that$$e^{-h s}=1-h s+\frac{(h s)^2}{2} e^{-h t}.$$
Thus
\allowdisplaybreaks
\begin{align*}
\|(\frac{e^{-hs}-1}{h}+s)e^{-\eta s}h(s) \| &\leq \frac{|h| s^2}{2}Me^{-(\eta-\xi) s}e^{-ht}\\
& \leq \frac{|h| s^2}{2}Me^{-(\eta-\xi) s}e^{|h|t}\\
&\leq\frac{|h| s^2}{2}Me^{-(\eta-\xi) s}e^{|h|s}\\
&=\frac{|h| s^2}{2}Me^{-(\eta-|h|-\xi) s}\\
& \leq \frac{|h| s^2}{2}Me^{-\gamma s}
\end{align*}
for any $s\geq 0$ and $h \in \bar{\bar{G}}$.
Consequently, for any $h \in \bar{\bar{G}}$, we have
\begin{align}
&\left\|\frac{H\left(\eta+h\right)-H\left(\eta\right)}{h}-\int_0^{\infty} e^{-\eta s}(-s)h(s)ds\right\|\nonumber\\
&= \left\|\int_0^{\infty}\left(\frac{e^{-hs}-1}{h}+s\right)e^{-\eta s}h(s)  d s\right\| \nonumber \\
&\leq  \int_0^{\infty}\left\|\left(\frac{e^{-hs}-1}{h}+s\right)e^{-\eta s}h(s) \right\|d s \nonumber \\
&\leq\frac{1}{2}|h| \int_0^{\infty}s^2Me^{-\gamma s}d s \nonumber \\
&= |h| \frac{M}{\gamma^3}. \nonumber
\end{align}
Letting $h \rightarrow 0$ in $\mathcal{T}_c$ in the above inequality, we have
$$\frac{H\left(\eta+h\right)-H\left(\eta\right)}{h}-\int_0^{\infty} e^{-\eta s}(-s)h(s)ds$$
converges  to  $\theta$ in $\mathcal{T}_c$, i.e.,
$$
\mathcal{T}_c\text{-}H^{\prime}(\eta)=\int_0^{\infty} e^{-\eta s}(-s)h(s)ds.
$$

\endproof

Based on Lemma 3.2, we can establish the Post-Widder inversion formula of such a Laplace  transform as follows.

\textbf{Theorem 3.1} (Post-Widder Inversion)
Suppose that $S$ is a $\mathcal{T}_{\varepsilon, \lambda}$-complete $R N$ module and $h:[0, \infty) \rightarrow (S, \mathcal{T}_{\varepsilon, \lambda})$ is a continuous function satisfying  $\|h(s)\| \leq M $ for any $s\geq 0$ and some $M \in L_{+}^0(\mathcal{F})$. For any $\eta \in L_{++}^0(\mathcal{F})$,
let $H(\eta)$ be the Laplace transform of $h$, then
$$
\lim _{k \rightarrow \infty}(-1)^k\left(\frac{k}{t}\right)^{k+1} \frac{1}{k!}\cdot\mathcal{T}_c\text{-}H^{(k)}\left(\frac{k}{t}\right)=h(t)
$$
in $ (S, \mathcal{T}_{\varepsilon, \lambda})$ for any $t>0$.

\textbf{Proof}
For any $\eta \in L_{++}^0(\mathcal{F})$, it is clear that $\eta \in  G_0$ and it follows from
Lemma 3.2 that
$$
\mathcal{T}_c\text{-}H^{(k)}(\eta)=\int_0^{\infty} e^{-\eta s}(-s)^k h(s)d s
$$
for any $k \in N$.  For any  $t>0$ and $k \in N$, we have $\frac{k}{t}>0$ and $\frac{k}{t}$ can be identified
with $\tilde{I}_{\Omega}\frac{k}{t}$ in $G_0$ conventionally,
thus
$$
\begin{aligned}
\mathcal{T}_c\text{-}H^{(k)}(\frac{k}{t})=&\int_0^{\infty} e^{-\frac{k}{t} s}(-s)^k h(s)d s\\
=&(-1)^{k}t^{k+1}\int_0^{\infty} e^{-k s} s^k h(ts)d s.
\end{aligned}
$$

Since $$\frac{k^{k+1}}{k !}\int_0^{\infty} e^{-ks} s^k  ds=1$$ for any $k \in N$, it follows that
$$
\begin{aligned}
J_k(t)&=(-1)^{k}(\frac{k}{t})^{k+1} \frac{1}{k !} \cdot\mathcal{T}_c\text{-}H^{(k)}(\frac{k}{t})-h(t)\\
& =\frac{1}{k !} k^{k+1} \int_0^{\infty} e^{-k s} s^k h(ts)d s-h(t)\\
& =\frac{1}{k !} k^{k+1} \int_0^{\infty} e^{-k s} s^k (h(ts)-h(t)) d s
\end{aligned}
$$
for any  $t>0$.

It is easy to see that
$$\|h(ts)-h(t)\|\leq  2M$$
for any $s, t\geq 0$. Now, let $\zeta=2M,$ then  $\zeta \in  L_+^0(\mathcal{F})$.
Set
$$ E_{n}=[n-1 \leq \zeta<n]$$
for any $n \in N$, then $E_{n}\in \mathcal{F}$, $E_{i}\bigcap E_{j}=\emptyset$ for any $i, j \in N$ with $i\neq j$, and further $\sum_{n=1}^{\infty} E_{n}=\Omega$.
Since, for any $n \in N$,  $\|I_{E_{n}}h(ts)-I_{E_{n}}h(t)\|\leq n$ for any $s, t\geq 0$ and $\|I_{E_{n}}h(ts)-I_{E_{n}}h(t)\|\rightarrow 0$ in probability $P$ as $s\rightarrow 1$, according to Lebesgue's dominated convergence theorem, we have
$\int_{\Omega}\|I_{E_{n}}h(ts)-I_{E_{n}}h(t)\|dP\rightarrow 0$  as $s\rightarrow 1$.
For any  $n \in N$ and $0<r<1$, due to  Proposition 2.3, we have
\allowdisplaybreaks[4] 
\begin{align}
&\frac{1}{k !} k^{k+1}\int_{\Omega}\left\|\int_{1-r}^{1+r} e^{-ks} s^k I_{E_{n}}(h(ts)-h(t))ds\right\|dP \nonumber \\
&\leq \frac{1}{k !} k^{k+1}\int_{\Omega} \int_{1-r}^{1+r}\left\|e^{-ks} s^k I_{E_{n}}(h(ts)-h(t))\right\|dsdP \nonumber \\
&= \frac{1}{k !} k^{k+1}\int_{1-r}^{1+r} \int_{\Omega}\left\|e^{-ks} s^k I_{E_{n}}(h(ts)-h(t))\right\| dPds \nonumber \\
&= \frac{1}{k !} k^{k+1}\int_{1-r}^{1+r}e^{-ks} s^k \int_{\Omega}\left\| I_{E_{n}}(h(ts)-h(t))\right\| dPds \nonumber \\
&\leq \max_{s\in [1-r, 1+r]}\int_{\Omega}\left\| I_{E_{n}}(h(ts)-h(t))\right\| dP\cdot\frac{k^{k+1}}{k !} \int_0^{\infty} e^{-ks} s^k ds \nonumber \\
&= \max_{s\in [1-r, 1+r]}\int_{\Omega}\left\| I_{E_{n}}(h(ts)-h(t))\right\| dP \nonumber \\
&\rightarrow 0\ \text{as}\ r\rightarrow 0, \nonumber
\end{align}
which implies that $ \frac{1}{k !} k^{k+1} \|\int_{1-r}^{1+r} e^{-ks} s^k I_{E_{n}}(h(ts)-h(t))ds\|\rightarrow 0$ in probability $P$ as $ r\rightarrow 0$.
Since  $\sum_{n=1}^{\infty} P(E_{n})=P\left(\sum_{n=1}^{\infty}E_{n}\right)=P(\Omega)=1$, it follows that $$ \frac{1}{k !} k^{k+1} \|\int_{1-r}^{1+r} e^{-ks} s^k(h(ts)-h(t))ds\|\rightarrow 0$$ in probability $P$ as $ r\rightarrow 0$, i.e.,
for any $\varepsilon, \lambda>0$, there exists a $\delta$ with $0<\delta<1$ such that
\begin{equation}
\begin{aligned}
P[\frac{1}{k !} k^{k+1} \|\int_{1-r}^{1+r} e^{-ks} s^k(h(ts)-h(t))ds\|\geq \frac{\varepsilon}{3}]\leq \frac{\lambda}{3}
\end{aligned}
\end{equation}
whenever $0<r\leq \delta$.

 Next, let $$J_k(t): =J_{k, 1}(t)+J_{k, 2}(t)+J_{k, 3}(t)$$
 for any $k \in N$ and $t>0$,
where the three integrals $J_{k, 1}(t), J_{k, 2}(t), J_{k, 3}(t)$ correspond to the intervals $(0, 1-\delta),(1-\delta, 1+\delta),(1+\delta, \infty)$, respectively. According to (1), for the above $\varepsilon$ and $\lambda$, we have
$$
\begin{aligned}
P[\left\|J_{k, 2}(t)\right\|\geq \frac{\varepsilon}{3}]\leq \frac{\lambda}{3}
\end{aligned}
$$
for any $k \in N$ and $t>0$.

It is clear that the mapping $s \rightarrow e^{-ks} s^{k}$ is increasing on $(0,1-\delta)$ for any $k \in N$.
Let $D=[M>0]$,  then we can assume that $P(D)>0$ without loss of generality. For any $k \in N$, we have
$$
\begin{aligned}
\left\|J_{k, 1}(t)\right\| &\leq \frac{k^{k+1}}{k!} e^{-k(1-\delta)}(1-\delta)^{k}\int_{0}^{1-\delta}\|h(ts)-h(t)\|ds\\
&\leq 2M\frac{k^{k+1}}{k!} e^{-k(1-\delta)}(1-\delta)^{k+1}\\
&= 2MI_{D}\frac{k^{k+1}}{k!} e^{-k(1-\delta)}(1-\delta)^{k+1}\\
&:=\xi_k.
\end{aligned}
$$
Since $\frac{\xi_{k+1}}{\xi_k}=e^{\delta}(1-\delta) \tilde{I}_{\Omega}<1\; on\;  \Omega$ for any $k \in N$, it follows that  $\|J_{k, 1}(t)\|\rightarrow 0$ in probability $P$ as $k\rightarrow \infty$ for any $t>0$, i.e.,
for the above $\varepsilon$ and $\lambda$, there exists an $N_1 \in N$ such that
 $$
\begin{aligned}
P\{\left\|J_{k, 1}(t)\right\|\geq \frac{\varepsilon}{3}]\leq \frac{\lambda}{3}
\end{aligned}
$$
as $k\geq N_1$.

Clearly,  the mapping $s \mapsto e^{-ks} s^k$ is decreasing on $(1+\delta, \infty)$ for any $k \in N$.
Furthermore, there exists an $N_2 \in N$ such that  for any $k > N_2$, we have

\allowdisplaybreaks[4] 
\begin{align*}
\left\|J_{k, 3}(t)\right\| & =\left\|\frac{k^{k+1}}{k !} \int_{1+\delta}^{\infty}  e^{-ks} s^k (h(ts)-h(t))d s\right\| \\
& \leq 2M\frac{k^{k+1}}{k !} \int_{1+\delta}^{\infty}e^{-ks}s^kd s\\
& \leq 2M\frac{k^{k+1}}{k !} e^{-(k-N_2)(1+\delta)}(1+\delta)^{k-N_2} \int_{1+\delta}^{\infty}e^{-N_2s}s^{N_2}d s\\
& = 2MI_{D}\frac{k^{k+1}}{k !} e^{-(k-N_2)(1+\delta)}(1+\delta)^{k-N_2} \int_{1+\delta}^{\infty}e^{-N_2s}s^{N_2}d s\\
& :=\zeta_{k}
\end{align*}
for any $t>0$.
Since $\frac{\zeta_{k+1}}{ \zeta_k} \rightarrow e^{-\delta}(1+\delta)\tilde{I}_{\Omega}<1\; on\;  \Omega$ as $k\rightarrow \infty$, it follows that
$\|J_{k, 3}(t)\|\rightarrow 0$ in probability $P$ as $k\rightarrow \infty$ for any $t>0$, i.e.,
for the above $\varepsilon$ and $\lambda$, there exists an $N_3 \in N$ with $N_3>N_2$ such that
 $$
\begin{aligned}
P[\left\|J_{k, 3}(t)\right\|\geq \frac{\varepsilon}{3}]\leq \frac{\lambda}{3}
\end{aligned}
$$
as $k\geq N_3$.

Consequently,  for any $t>0$ and the above $\varepsilon$, $\lambda$, one can obtain
$$
\begin{aligned}
P[\left\|J_k(t)\right\|\geq \varepsilon]&\leq P[\left\|J_{k, 1}(t)\right\|\geq \frac{\varepsilon}{3}]+P[\left\|J_{k, 2}(t)\right\|\geq \frac{\varepsilon}{3}]+P[\left\|J_{k, 3}(t)\right\|\geq \frac{\varepsilon}{3}]\\
&\leq \lambda
\end{aligned}
$$
as $k\geq \max\{N_1, N_3\}$, which shows that
$$
\lim _{k \rightarrow \infty}(-1)^k\left(\frac{k}{t}\right)^{k+1} \frac{1}{k!}\cdot\mathcal{T}_c\text{-}H^{(k)}\left(\frac{k}{t}\right)=h(t)
$$
in $(S, \mathcal{T}_{\varepsilon, \lambda})$.
\endproof

\textbf{Corollary 3.1} (Uniqueness)
Suppose that $S$ is a $\mathcal{T}_{\varepsilon, \lambda}$-complete $R N$ module and $h_1, h_2:[0, \infty) \rightarrow (S, \mathcal{T}_{\varepsilon, \lambda})$ are two continuous functions satisfying  $\|h_1(s)\| \leq M_1 e^{\xi_1 s}$ and $\|h_2(s)\| \leq M_2 e^{\xi_2 s}$ for any $s\geq 0$  and some $M_1, M_2 \in L_{+}^0(\mathcal{F}), \xi_1, \xi_2 \in L^0(\mathcal{F}, R)$. If
$$\int_0^{\infty} e^{-\eta s}h_1(s)ds=\int_0^{\infty} e^{-\eta s}h_2(s)ds$$
for any $\eta \in L^0(\mathcal{F}, R)$ with $\eta>\bigvee \{\xi_1, \xi_2\}$ on  $\Omega$, then $h_1=h_2$.

\textbf{Proof}
Let $h(s)=e^{-\bigvee \{\xi_1, \xi_2\} s}(h_1(s)-h_2(s))$ for any $s\geq 0$, then $\|h(s)\|\leq M_1+M_2$.  Thus the Laplace transform of $h$ can be given by
$$
\begin{aligned}
H(\zeta)  =\int_0^{\infty} e^{-\zeta s} h(s) d s \\
\end{aligned}
$$
for any $\zeta\in G_0$. Since

$$\int_0^{\infty} e^{-\eta s}(h_1(s)-h_2(s))ds=\theta$$
for any $\eta \in L^0(\mathcal{F}, R)$ with $\eta>\bigvee \{\xi_1, \xi_2\}$ on  $\Omega$, it follows that

\allowdisplaybreaks
\begin{align*}
H(\zeta) & =\int_0^{\infty} e^{-\zeta s} h(s) d s \\
& =\int_0^{\infty} e^{-\zeta s} e^{-\bigvee \{\xi_1, \xi_2\} s}\left(h_1(s)-h_2(s)\right) d s \\
& =\int_0^{\infty} e^{-(\zeta+\bigvee \{\xi_1, \xi_2\}) s}\left(h_1(s)-h_2(s)\right) d s \\
& =\theta
\end{align*}
for any $\zeta\in G_0 $.
Since $\|h(s)\|\leq M_1+M_2$ for any $s\geq 0$, it follows from Theorem 3.1 that $h(s)=e^{-\bigvee \{\xi_1, \xi_2\} s}(h_1(s)-h_2(s))=\theta$, i.e., $h_1(s)-h_2(s)=\theta$. Thus $h_1=h_2$.
\endproof

\section{$C$-semigroups of continuous module homomorphisms on a  $\mathcal{T}_{\varepsilon, \lambda}$-complete $R N$ module}

\quad \, In the sequel of this paper, we always assume that $(S,\|\cdot\|)$ is a  $\mathcal{T}_{\varepsilon, \lambda}$-complete $R N$ module over K with base $(\Omega, \mathcal{F}, P)$. The main result of this section is Theorem 4.1, in which we establish the Hille-Yosida type
theorem for exponentially bounded $C$-semigroups on a  $\mathcal{T}_{\varepsilon, \lambda}$-complete $R N$ module for the case that $R(C)$
is  nondense in  $(S, \mathcal{T}_{\varepsilon, \lambda})$.

\textbf{Definition 4.1} \textrm{(see \cite{Son})}
Let $C\in B(S)$ be an injective operator on $S$. Then a  family $\{W(t): t \geq 0\} \subset B(S)$ is called a $C$-semigroup on $S$ if

(1) $CW(s+t)=W(t) W(s)$ for any $t, s\geq 0$;

(2) $W(0)=C$;

(3) $W(s)$ is strongly continuous, i.e., for any $x \in S$, the mapping $s \rightarrow W(s) x$ from $[0, \infty)$ into $(S, \mathcal{T}_{\varepsilon, \lambda})$ is continuous.

\textbf{Definition 4.2} \textrm{(see \cite{Zhangliuwei})}
 Suppose that $\{W(s):s \geq 0\}$ is a $C$-semigroup on $S$ and $R(C)$ denotes the range of $C$.  Define
$$
D(A)=\left\{x \in S: C^{-1} \lim _{s \rightarrow 0} \frac{W(s) x-C x}{s} \text { exists and belongs to } R(C)\right\}
$$
and
$$
Ax=C^{-1} \lim _{s \rightarrow 0} \frac{W(s) x-C x}{s}
$$
for any $x \in D(A)$, then the mapping $A: D(A) \rightarrow S$ is called the infinitesimal generator of $\{ W(s): s \geq 0\}$, also denoted by $(A, D(A))$ in this section.

In this paper, a $C$-semigroup $\{W(s): s \geq 0\}$ is said to be locally almost surely bounded if for any $l>0, \bigvee_{s \in[0, l]}\|W(s)\|$ is in $L_{+}^0(\mathcal{F})$. Besides, $\{W(s): s \geq 0\}$ is said to be exponentially bounded on $S$ if there are $M \in L_{+}^0(\mathcal{F}), \xi \in L^0(\mathcal{F}, R)$ satisfying
$\|W(s)\| \leq M e^{\xi s}$
for any $s \geq 0$, also denoted by $H(M, \xi)$ the set of all such exponentially bounded $C$-semigroups $\{W(s):s \geq 0\}$ on $S$ in this paper.

\textbf{Proposition 4.1} \textrm{(see \cite{Zhangliuwei})}
Let $\{W(s): s \geq 0\}$ be a locally almost surely bounded $C$-semigroup on $S$ with the infinitesimal generator $(A, D(A))$. Then

(1) for any $x \in D(A)$ and  $s \geq 0$, we have $W(s) x \in D(A)$ and
$$
\frac{dW(s) x}{ds}=W(s) Ax=A W(s) x;
$$

(2) for any $x \in S$ and $s \geq 0$, we have $\int_0^s W(t) x dt \in D(A)$ and
$$
A \int_0^s W(t) x dt=W(s)x-Cx;
$$

(3) $R(C) \subseteq \overline{D(A)}^{\mathcal{T}_{\varepsilon, \lambda}}$;

(4) $A$ is closed and $C^{-1} A C=A$;

(5) for any $x \in$ $D(A)$, the mapping $s \rightarrow C^2 W(s)x$ is locally $L^0$-Lipschitz.

Besides, just as in the classical case, Proposition 4.2 below shows that the infinitesimal generator of a locally almost surely bounded $C$-semigroup determines the $C$-semigroup uniquely.

\textbf{Proposition 4.2}
Let $\{W(s):s \geq 0\}$ and  $\{S(s):s \geq 0\}$ be  two locally almost surely bounded $C$-semigroups with the infinitesimal generators $(A, D(A))$ and $(B, D(B))$,  respectively. If $A=B$, then $W(s)=S(s)$ for any $ s\geq 0$.

\textbf{Proof}
For any $x \in S$ and $t \geq 0$, define a mapping $g: [0, t] \rightarrow S$ by $$g(s)=C^2W(t-s) \int_0^s S(r) x d r,$$ then it follows from Propositions 2.2 and 4.1(2) that
$$
\begin{aligned}
g^{\prime}(s) & =\frac{d}{d s}\left[C^2W(t-s) \int_0^s S(r) x d r\right] \\
& =-C^2W(t-s) A \int_0^s S(r) x d r+C^2W(t-s)S(s)x \\
& =-C^2W(t-s)(S(s)x-Cx)+C^2W(t-s)S(s)x \\
& =C^3W(t-s)x
\end{aligned}
$$
for any $s \in [0, t]$.
Moreover, for any $x \in S$ and $t\geq 0$, we have
\allowdisplaybreaks[4] 
\begin{align}
\|g(s_1)-g(s_2)\|\nonumber
&= \left\| C^2 W(t-s_1) \int_0^{s_1} S(r)x \, dr - C^2 W(t-s_2) \int_0^{s_2} S(r)x \, dr \right\| \nonumber \\
&= \left\| C^2 W(t-s_1) \int_0^{s_1} S(r)x \, dr - C^2 W(t-s_2) \int_0^{s_1} S(r)x \, dr \right. \nonumber \\
&\quad \left. - C^2 W(t-s_2) \int_{s_1}^{s_2} S(r)x \, dr \right\| \nonumber \\
&\leq \left\| C^2 W(t-s_1) \int_0^{s_1} S(r)x \, dr - C^2 W(t-s_2) \int_0^{s_1} S(r)x \, dr \right\| \nonumber \\
&\quad + \left\| C^2 W(t-s_2) \int_{s_1}^{s_2} S(r)x \, dr \right\| \nonumber \\
&\leq \left\| C^2 W(t-s_1) \int_0^{s_1} S(r)x \, dr - C^2 W(t-s_2) \int_0^{s_1} S(r)x \, dr \right\| \nonumber \\
&\quad + \bigvee_{u \in [0,t]} \|C^2 W(u)\| \bigvee_{s \in [s_1,s_2]} \|S(s)x\| |s_1-s_2| \nonumber
\end{align}
for any $s_1, s_2 \in [0,t]$.
 According to Proposition 4.1(5), it follows that $g$ is locally $L^0$-Lipschitz.
Further, due to Proposition 2.1, we have
$$
C^3\int_0^t S(r) x d r-0=\int_0^tg^{\prime}(s)ds=C^3\int_0^tW(t-s)xds
 $$
for any $x \in S$ and $t\geq 0$, i.e.,
$$
C^3\int_0^t W(r) x d r=C^3\int_0^t S(r) x d r.
 $$
Since $C$ is injective, it follows that $W(s)= S(s)$ for any $s \geq 0$.
\endproof

It should be pointed out that an exponentially bounded $C$-semigroup in a $\mathcal{T}_{\varepsilon, \lambda}$-complete $R N$ module is locally almost surely bounded. However, due to  Example 3.3  in \cite{Son}, a locally almost surely bounded $C$-semigroup in a $\mathcal{T}_{\varepsilon, \lambda}$-complete $R N$ module may not be exponentially bounded.

Subsequently, we always assume that $C$ belongs to $B(S)$ and is an injective operator on $S$.
Suppose that $(A, D(A))$ is a module homomorphism on $S$ and  $\rho_C(A)$ denotes the set $\{\eta \in L^0(\mathcal{F}, K):\eta-A$ is injective and $\left.R(C) \subseteq R(\eta-A)\right\}$, where $R(\eta-A)$ stands for the range of the module homomorphism  $(\eta-A)$,  then $\rho_C(A)$ is called the random $C$-resolvent set of $A$. Moreover, if $\eta \in \rho_C(A)$, then $R_C(\eta, A):=(\eta-A)^{-1} C$ is called the random $C$-resolvent of $A$.

\textbf{Lemma 4.1}
Suppose that $\{W(s):s \geq 0\}$ belongs to $ H(M, \xi)$ and $(A, D(A))$ is the infinitesimal generator of $\{W(s):s \geq 0\}$.
If $\eta \in G_\xi$, then $\eta \in \rho_C(A)$ and
$$R_C(\eta, A)x=\int_0^\infty e^{-\eta s}W(s)xds$$
for any $x \in S$.

\textbf{Proof}
For the readers' convenience, the proof is divided into two steps.

\vspace{1\baselineskip}
Step 1.
For any  $\eta \in G_\xi$ and $x_0 \in S$, we have
$$
\begin{aligned}
\|e^{-\eta s} W(s) x_0-Cx_0\|&\leq  \|e^{-\eta s} W(s) x_0\|+ \|Cx_0\|\\
&\leq M e^{-(\eta-\xi) s}\|x_0\|+\|Cx_0\|\\
&\leq M\|x_0\|+\|Cx_0\|\\
\end{aligned}
$$
for any $s\geq 0$. Now, let
$$
\zeta_{x_0}=M\|x_0\|+\|Cx_0\|,
$$
then $\zeta_{x_0} \in  L_+^0(\mathcal{F})$.
Set
$$ E_{k, x_0}=[k-1 \leq \zeta_{x_0}<k]$$
for any $k \in N$, then $E_{k, x_0}\in \mathcal{F}$, $E_{i, x_0}\bigcap E_{j, x_0}=\emptyset$ for any $i, j \in N$ with $i\neq j$, and further $\sum_{k=1}^{\infty} E_{k, x_0}=\Omega$.
Clearly, for any $k \in N$ and $\eta \in G_\xi$, $\|I_{E_{k, x_0}}e^{-\eta s}W(s) x_0-I_{E_{k, x_0}}Cx_0\|\leq k$ for any $s\geq 0$, and further $$\|I_{E_{k, x_0}}e^{-\eta s}W(s) x_0-I_{E_{k, x_0}}Cx_0\|\rightarrow 0$$ in probability $P$ as $s\rightarrow 0$. Thus, according to Lebesgue's dominated convergence theorem, we have
$$\int_{\Omega}\|I_{E_{k, x_0}}e^{-\eta s}W(s) x_0-I_{E_{k, x_0}}Cx_0\|dP\rightarrow 0$$  as $s\rightarrow 0$.
Consequently,  due to  Proposition 2.3, we have
\allowdisplaybreaks[4] 
\begin{align}
& \frac{1}{l}\int_{\Omega}\left\|\int_0^l I_{E_{k, x_0}}e^{-\eta s}W(s) x_0-I_{E_{k,x_0}}Cx_0ds\right\|dP \nonumber \\
&\leq  \frac{1}{l}\int_{\Omega} \int_0^l\left\|I_{E_{k, x_0}}e^{-\eta s}W(s) x_0-I_{E_{k,x_0}}Cx_0\right\|dsdP \nonumber \\
&=  \frac{1}{l}\int_0^l \int_{\Omega}\left\|I_{E_{k, x_0}}e^{-\eta s}W(s) x_0-I_{E_{k,x_0}}Cx_0\right\| dPds \nonumber \\
&\leq \max_{s\in [0, l]}\int_{\Omega}\left\|I_{E_{k, x_0}}e^{-\eta s}W(s) x_0-I_{E_{k,x_0}}Cx_0\right\| dP \nonumber \\
&\rightarrow 0\ \text{as}\ l\rightarrow 0, \nonumber
\end{align}
which implies that $\frac{1}{l} \|\int_0^l(I_{E_{k, x_0}}e^{-\eta s}W(s)x_0-I_{E_{k, x_0}}Cx_0)ds\|\rightarrow 0$ in probability $P$ as $l\rightarrow 0$.
Since $$\sum_{k=1}^{\infty} P(E_{k,x_0})=P(\sum_{k=1}^{\infty}E_{k,x_0})=P(\Omega)=1,$$ it follows that
\begin{equation}
\begin{aligned}
\lim _{l\rightarrow 0} \frac{1}{l}\int_0^l e^{-\eta s}  W(s) x_0 d s=C x_0
\end{aligned}
\end{equation}
for any  $\eta \in G_\xi$ and $x_0 \in S$.

Step 2.
For any $\eta \in G_\xi$ and $l>0$, we have
$$
\begin{aligned}
& \frac{W(l) \int_0^{\infty} e^{-\eta s} W(s) x d s-C\int_0^{\infty} e^{-\eta s} W(s) x d s}{l} \\
& =\frac{1}{l} \int_0^{\infty} e^{-\eta s}C W(s+l) x ds-\frac{1}{l}C \int_0^{\infty} e^{-\eta s} W(s) x ds \\
& =\frac{1}{l} C\int_l^{\infty} e^{-\eta(s-l)}W(s) x ds-\frac{1}{l} C\int_0^{\infty} e^{-\eta s}W(s) x d s \\
& =\frac{1}{l}\left(e^{\eta l}-1\right) C\int_0^{\infty} e^{-\eta s} W(s) x ds-\frac{1}{l} e^{\eta l} C\int_0^l e^{-\eta s} W(s) x ds
\end{aligned}
$$
for any $x \in S$.
Letting $l\rightarrow 0$ in the above equality, due to (2),
one can obtain
$$
CA \int_0^{\infty} e^{-\eta s} W(s) x ds=C\eta \int_0^{\infty} e^{-\eta s} W(s) x ds-C^2x.
$$
Since $C$ is injective, it follows that
$$
A \int_0^{\infty} e^{-\eta s} W(s) x d s=\eta \int_0^{\infty} e^{-\eta s} W(s) x d s-C x.
$$
Thus
\begin{equation}
\begin{aligned}
(\eta-A)\left[\int_0^{\infty} e^{-\eta s} W(s) x d s\right]=C x
\end{aligned}
\end{equation}
for any $\eta \in G_\xi$, which implies that $R(C) \subseteq R(\eta-A)$.

\vspace{1\baselineskip}
According to Proposition 4.1, we have $A$ is closed and
$$
\frac{dW(s) x}{ds}=W(s) Ax=A W(s) x
$$
for any $x \in D(A)$, thus
$$
A \int_0^{\infty} e^{-\eta s} W(s) x d s=\int_0^{\infty} e^{-\eta s} W(s) A x d s
$$
for any $\eta \in G_{\xi}$ and $x \in D(A)$.
Consequently, due to (3), we have
\begin{equation}
\begin{aligned}
\int_0^{\infty} e^{-\eta s} W(s)(\eta-A) x d s&=(\eta-A)\left[\int_0^{\infty} e^{-\eta s} W(s) x d s\right] \\
&=C x
\end{aligned}
\end{equation}
for any  $x \in D(A)$. Moreover, for any $\eta \in G_{\xi}$, if $(\eta-A) x=\theta$ for some $x \in D(A)$, then

$$
C x=\int_0^{\infty} e^{-\eta s} W(s)(\eta-A) x d s=0.
$$
Hence $x=\theta$ since $C$ is injective, which shows that $\eta-A$ is injective for any $\eta \in G_{\xi}$. According to (3), one can obtain
$$
\begin{aligned}
&(\eta-A)\left[\int_0^{\infty} e^{-\eta s} W(s) x d s-(\eta-A)^{-1} C x\right]\\
&=(\eta-A)\left[\int_0^{\infty} e^{-\eta s} W(s) x d s-R_C(\eta, A) x\right]\\
&=0.
\end{aligned}
$$
For any $\eta \in G_{\xi}$, since $\eta-A$ is injective and $R(C) \subseteq R(\eta-A)$, it follows that
$\eta \in \rho_C(A)$ and
$$
R_C(\eta, A) x=\int_0^{\infty} e^{-\eta s} W(s) x d s
$$
for any $x \in S$, which completes the proof of Lemma 4.1.

\endproof

Based on Lemma 4.1 and the uniqueness of the Laplace transform established in Corollary 3.1, we can now prove Theorem 4.1 below, which establishes the Hille-Yosida type theorem for the case when $R(C)$ is nondense in $(S, \mathcal{T}_{\varepsilon, \lambda})$.

\textbf{Theorem 4.1} 
Suppose that $(A, D(A))$ is a module homomorphism on $S$ and $W:[0, \infty) \rightarrow B(S)$ is a strongly continuous family satisfying $\|W(s)\| \leq M e^{\xi s}$ for any $s\geq 0$ and some $M \in L_{+}^0(\mathcal{F})$, $\xi \in L^0(\mathcal{F}, R)$. Then the following statements are equivalent.

$(a)$ $\{W(s): s \geq 0\}$ is a $C$-semigroup on $S$ with the infinitesimal generator $(A, D(A))$.

$(b)$ $(b_1)$ $A=C^{-1} A C$;

~~~~$(b_2)$ If $\eta \in G_{\xi}$, then $\eta \in \rho_C(A)$ and
$$
R_C(\eta, A) x=\int_0^{\infty} e^{-\eta s} W(s)x ds
$$
for any $x \in S$.

\textbf{Proof}
$(\mathrm{a}) \Rightarrow(\mathrm{b})$.
According to Proposition 4.1(4) and Lemma 4.1, it follows that $(b)$ holds.

(b) $\Rightarrow(\mathrm{a})$. For any $\eta \in G_{\xi}$ and $x \in S$, define two functions $h_1, h_2: [0, \infty) \rightarrow S$
by $$h_1(t) = \int_0^{\infty} e^{-\eta s} W(t) W(s) x  ds$$ and $$h_2(t) = \int_0^{\infty} e^{-\eta s} W(t+s) Cx  ds,$$
then $$h_1(t) = W(t) \int_0^{\infty} e^{-\eta s} W(s) x  ds$$ since $W(t) \in B(S)$ for any  $t \geq 0$.
Further,  $h_1: [0, \infty) \rightarrow (S, \mathcal{T}_{\varepsilon, \lambda})$ is continuous and
\allowdisplaybreaks
\begin{align*}
\|h_1(t)\| & = \| W(t) \int_0^{\infty} e^{-\eta s} W(s) x  ds \| \\
& \leq M e^{\xi t} \| \int_0^{\infty} e^{-\eta s} W(s) x  ds \|\\
& \leq \frac{M^2\|x\|}{\eta-\xi} e^{\xi t}.
\end{align*}
For $h_2$, we have
\allowdisplaybreaks[4] 
\begin{align}
\|h_2(t)\| & = \| \int_0^{\infty} e^{-\eta s} W(t+s) Cx  ds \| \nonumber \\
& \leq \int_0^{\infty} e^{-\eta s} \|W(t+s) Cx\|  ds \nonumber \\
& \leq M e^{\xi t} \|Cx\| \int_0^{\infty} e^{-(\eta - \xi) s}  ds \nonumber \\
&=\frac{M\|Cx\|}{\eta-\xi}e^{\xi t}. \nonumber
\end{align}
Next, we will prove that $h_2: [0, \infty) \rightarrow (S, \mathcal{T}_{\varepsilon, \lambda})$ is continuous.
Clearly,
\[
\begin{aligned}
h_2(t)&=\int_0^{\infty} e^{-\eta s} W(t+s) Cx d s\\
&=e^{\eta t} \int_t^{\infty} e^{-\eta v} W(v)C x d v.
\end{aligned}
\]
Further, using the inequality $$|e^{\zeta_{1}} - e^{\zeta_{2}}| \leq \frac{1}{2}|\zeta_{1}-\zeta_{2}|(e^{\zeta_{1}} + e^{\zeta_{2}})$$ for any $\zeta_{1}, \zeta_{2} \in L^0(\mathcal{F}, R)$, one can obtain
\allowdisplaybreaks[4] 
\begin{align}
&\|h_2(t) - h_2(t_0)\| \nonumber \\
& =\| e^{\eta t} \int_t^{\infty} e^{-\eta v} W(v) C x  dv - e^{\eta t_0} \int_{t_0}^{\infty} e^{-\eta v} W(v) C x  dv \| \nonumber \\
& \leq e^{\eta t} \| \int_t^{\infty} e^{-\eta v} W(v) C x  dv - \int_{t_0}^{\infty} e^{-\eta v} W(v) C x  dv \| \nonumber \\
&\quad + \| \int_{t_0}^{\infty} e^{-\eta v} W(v) C x  dv \| \cdot |e^{\eta t} - e^{\eta t_0}| \nonumber \\
& = e^{\eta t} \| \int_t^{t_0} e^{-\eta v} W(v) C x  dv \| \nonumber \\
&\quad + \| \int_{t_0}^{\infty} e^{-\eta v} W(v) C x  dv \| \cdot |e^{\eta t} - e^{\eta t_0}| \nonumber \\
& \leq e^{\eta t}M \|C x\| \cdot \frac{| e^{-(\eta-\xi) t} - e^{-(\eta-\xi) t_0} |}{\eta - \xi} \nonumber \\
&\quad +  \int_{t_0}^{\infty} e^{-\eta v} \|W(v) C x \| dv  \cdot |e^{\eta t} - e^{\eta t_0}| \nonumber \\
& \leq \frac{1}{2} M \|C x\| |t - t_0| e^{\eta t} \left( e^{-(\eta-\xi) t} + e^{-(\eta-\xi) t_0} \right) \nonumber \\
&\quad + \frac{1}{2} \frac{M\eta}{\eta-\xi}\|C x\| |t - t_0| e^{-(\eta-\xi) t_0}(e^{\eta t} + e^{\eta t_0}) \nonumber\\
& \leq  M \|C x\| |t - t_0| e^{\eta t}  + \frac{1}{2} \frac{M\eta}{\eta-\xi}\|C x\| |t - t_0|(e^{\eta t} + e^{\eta t_0}) \nonumber
\end{align}
for any $t, t_0 \geq 0$. In particular, if we choose $t\in [\frac{t_0}{2}, \frac{3t_0}{2} ]$, then
\[
\begin{aligned}
\|h_2(t) - h_2(t_0)\|&\leq  M \|C x\| |t - t_0| e^{\frac{3}{2}\eta t_0}\\
 &\quad + \frac{1}{2} \frac{M\eta}{\eta-\xi}\|C x\| |t - t_0| (e^{\frac{3}{2}\eta t_0} + e^{\eta t_0})
 \end{aligned}
\]
which implies that $\|h_2(t) - h_2(t_0)\| \rightarrow 0$ in probability $P$ as $t \rightarrow t_0$,
i.e., $h_2: [0, \infty) \rightarrow (S, \mathcal{T}_{\varepsilon, \lambda})$ is continuous.

Since $R_C(\eta, A)C-R_C(\mu, A)C=(\mu-\eta) R_C(\mu, A) R_C(\eta, A)$ for any  $\mu,\ \eta \in G_{\xi}$, it follows that
\allowdisplaybreaks[4] 
\begin{align}
    \int_0^{\infty} e^{-\mu t}h_2(t)  dt &= \int_0^{\infty} e^{-\mu t}\int_0^{\infty}  e^{-\eta s} W(t+s) C x  ds  dt \nonumber \\
    &= \int_0^{\infty} \int_0^{\infty} e^{-\mu t} e^{-\eta s} W(t+s) C x  ds  dt \nonumber \\
    &= \int_0^{\infty} e^{-(\mu-\eta) t}\left(R_C(\eta, A) C x-\int_0^t e^{-\eta v}W(v) C x  d v\right)  d t \nonumber \\
    &= \frac{1}{\mu-\eta}\left(R_C(\eta, A) C x-R_C(\mu, A) C x\right) \nonumber \\
    &= R_C(\mu, A) R_C(\eta, A) x \nonumber \\
    &= \int_0^{\infty} \int_0^{\infty} e^{-\mu t} e^{-\eta s} W(t)W(s) x  d s  d t \nonumber \\
    &= \int_0^{\infty} e^{-\mu t}\int_0^{\infty}  e^{-\eta s} W(t)W(s) x  d s  d t \nonumber \\
    &= \int_0^{\infty} e^{-\mu t}h_1(t)  dt \nonumber
\end{align}
for any  $x \in S$ and $ \mu>\eta$ on  $\Omega$.
According to Corollary 3.1, we have $h_1=h_2$, i.e.,
 $$\int_0^{\infty} e^{-\eta s} W(t+s) Cx  ds=\int_0^{\infty}  e^{-\eta s} W(t)W(s) x  d s.$$
Applying Corollary 3.1 again, it follows that $$W(t+s) C=W(t)W(s)$$ for any $t, s \geq 0$.
If $W(0) x=\theta$ for some $x \in S$, then $$W(s) C x=W(s) W(0) x=\theta.$$ Further, according to $(b_2)$, we have
\begin{align}\nonumber
(\eta-A)^{-1} C^2x&=R_C(\eta, A) C x\\\nonumber
&=\int_0^{\infty} e^{-\eta s} W(s)Cx ds\\\nonumber
&=\theta\nonumber
\end{align}
for any $\eta \in G_{\xi}$.  Thus $$(\eta-A)(\eta-A)^{-1} C^2x=(\eta-A)\theta,$$
i.e., $C^2x=\theta.$ Since $C$ is injective, it follows that $x=\theta$, i.e.,  $W(0)$ is  injective.
Due to $$W(0)(W(0)-C)=0,$$ we have $W(0)=C$.
 Thus $\{W(s):s \geq 0\}$ is a $C$-semigroup on $S$.

 Finally, we will show that  $(A, D(A))$ is the infinitesimal generator of $\{W(s): s\geq 0\}$. Suppose that $(B, D(B))$
is the infinitesimal generator of $\{W(s): s\geq 0\}$. Due to Proposition 4.1(4), we obtain
$C^{-1}BC=B$. For any $x \in D\left(C^{-1} A C\right)$, let $$y=\eta x-C^{-1} A C x$$ for any $\eta \in G_{\xi}$, then
\[
\begin{aligned}
C x=(\eta- A)^{-1}Cy=R_C(\eta, A) y.
\end{aligned}
\]
 According to Lemmas 4.1 and $(b_2)$, we have
\[
\begin{aligned}
R_C(\eta, B) y=\int_0^{\infty} e^{-\eta s} W(s)yds=R_C(\eta, A) y
\end{aligned}
\]
for any $\eta \in G_{\xi}$.
Thus
\[
\begin{aligned}
C x=R_C(\eta, B) y=(\eta- B)^{-1}Cy
\end{aligned}
\]
for any $\eta \in G_{\xi}$,
which shows that $$y=\eta x-C^{-1} B C x.$$ Consequently, we can obtain $$C^{-1} A C \subseteq C^{-1}B C.$$
Similarly, one can prove that $C^{-1} B C \subseteq C^{-1} A C$. Then $C^{-1} B C= C^{-1} A C$.
Since $A=C^{-1} A C$ and $C^{-1}BC=B$, it follows that $A=B$.
\endproof

If we choose $\mathcal{F}=\{\Omega, \emptyset\}$, then the $\mathcal{T}_{\varepsilon, \lambda}$-complete $R N$ module $S$ reduces to a Banach space $X$ and the exponentially bounded $C$-semigroup $\{W(s): s \geq 0\}$ on $S$ reduces to an ordinary exponentially bounded $C$-semigroup on $X$, which leads to the following Corollary
4.1.

\textbf{Corollary 4.1} \textrm{(see \cite{deLaubenfels1})}
Let $(A, D(A))$ be a linear operator on  a Banach space $X$, $B(X)$ the space of all bounded linear operators on $X$, and $W: [0, \infty) \rightarrow B(X)$ a strongly continuous family satisfying $\|W(s)\| \leq M e^{a s}$ for any $s \geq 0$ and some $M \geq 0, a \in R$. Then the following statements are equivalent.

$(a)$ $\{W(s): s \geq 0\}$ is a $C$-semigroup on $X$ with the infinitesimal generator $(A, D(A))$.

$(b)$ $(b_1)$ $A=C^{-1} A C$;

~~~~$(b_2)$ If $\eta >a$, then $\eta \in \rho_C(A)$ and
$$
R_C(\eta, A) x=\int_0^{\infty} e^{-\eta t} W(t) x d t
$$
for any $x \in X$, where $\rho_C(A)$ denotes the $C$-resolvent set of $A$ and $R_C(\eta, A)$ denotes the $C$-resolvent of $A$.

\section{An application to abstract Cauchy problems on a $\mathcal{T}_{\varepsilon, \lambda}$-complete $RN$ module}

~\quad In this section, let $A$ be a module homomorphism from $D(A)$ into $S$ and $[D(A)]$ denote the $RN$ module $D(A)$ with the following graph $L^0$-norm
$$
\|x\|_{[D(A)]} :=\|x\|+\|Ax\|
$$
for any $x \in D(A)$. Besides, $C([0, \infty),[D(A)])$ denotes the set of continuous  functions from $[0, \infty)$ to $[D(A)], C([0, \infty), S)$ the set of continuous  functions from $[0, \infty)$ to $S$, and $C^1([0, \infty), S)$ the set of continuously differentiable  functions from $[0, \infty)$ to $S$.

Suppose that $\{W(s):s \geq 0\}$ belongs to $ H(M, \xi)$ and $(A, D(A))$ is the infinitesimal generator of $\{W(s):s \geq 0\}$. In Lemma 4.1,
the random $C$-resolvent $R_C(\eta, A)$ of its infinitesimal generator $(A, D(A))$ can be characterized by the Laplace transform of the $C$-semigroup as follows
\begin{equation}
R_C(\eta, A)x = \int_0^{\infty} e^{-\eta s} W(s)x ds
\end{equation}
for any  $x \in S$ and  $\eta \in G_{\xi}$. In this section, we will apply (5) to Theorem 5.1.

\textbf{Theorem 5.1}
Suppose that $\{W(s):s \geq 0\}$ belongs to $ H(M, \xi)$ and $(A, D(A))$ is the infinitesimal generator of $\{W(s):s \geq 0\}$. Then, under the
locally $L^0$-Lipschitz condition on the solution, the abstract Cauchy problem
\begin{equation}
\begin{aligned}
\left\{\begin{array}{l}
\frac{du(t)}{d t}=Au(t),  \forall t>0, \\
u(0)=u_0\in R(R_C(\eta, A))
\end{array}\right.
\end{aligned}
\end{equation}
has a unique solution $u(t):=W(t)C^{-1}u_0$ belonging to  $C([0, \infty),[D(A)]) \cap C^1([0, \infty), S)$, where $\eta \in G_{\xi}$ and $R\left(R_C(\eta, A)\right)$ denotes the range of $R_C(\eta, A)$.

\textbf{Proof}
For any  $\eta \in G_{\xi}$ and $u_0\in R(R_C(\eta, A))$, there exists a $y_0\in S$ such that $u_0= R_C(\eta, A)y_0=(\eta-A)^{-1}Cy_0$.
According to (5), we have
$$
\begin{aligned}
u_0 & =\int_0^{\infty} e^{-\eta s}W(s)y_0 d s.
\end{aligned}
$$
Let $u(t)=W(t) C^{-1} u_0$ for any $t \geq 0$, then
\allowdisplaybreaks[4] 
\begin{align*}
u(t)&=W(t)C^{-1}\int_0^{\infty} e^{-\eta s}W(s)y_0 d s\\
&=\int_0^{\infty} e^{-\eta s}W(s+t)y_0 d s\\
&=e^{\eta t} \int_t^{\infty}e^{-\eta s}W(s)y_0 d s.
\end{align*}
Thus
$$u^{\prime}(t)=\eta \cdot e^{\eta t} \int_t^{\infty} e^{-\eta s} W(s)y_0 d s-W(t)y_0$$
for any $t \geq 0$.

For any $t \geq 0$, we have
\allowdisplaybreaks[4] 
\begin{align}
& \frac{W(l)e^{\eta t} \int_t^{\infty} e^{-\eta s}W(s)y_0  ds-Ce^{\eta t} \int_t^{\infty} e^{-\eta s}W(s)y_0  ds}{l} \nonumber \\
& = \frac{1}{l}\left[Ce^{\eta t} \int_t^{\infty} e^{-\eta s} W(l+s) y_0  ds-Ce^{\eta t} \int_t^{\infty} e^{-\eta s}W(s)y_0  ds\right] \nonumber \\
& =  \frac{1}{l}\left[Ce^{\eta t} \int_{t+l}^{\infty} e^{-\eta (s-l)} W(s)y_0  ds-Ce^{\eta t} \int_t^{\infty} e^{-\eta s} W(s)y_0  ds\right] \displaybreak[2] \nonumber \\ 
& =  \frac{1}{l}C \cdot e^{\eta t}\left[\int_t^{\infty}\left(e^{-\eta(s-l)}-e^{-\eta s}\right)W(s)y_0  ds-\int_t^{t+l} e^{-\eta(s-l)}W(s) y_0  ds\right] \nonumber \\
& = \frac{1}{l} e^{\eta t}\left(e^{\eta l}-1\right)C \int_t^{\infty}e^{-\eta s}W(s)y_0  ds-\frac{1}{l}W(t) \cdot e^{\eta l} \int_0^l e^{-\eta s}W(s)y_0  ds \nonumber
\end{align}
for any $l>0$, letting $l \rightarrow 0$  in the above equality, according to (2), we have
$$
CAu(t)=\eta\cdot e^{\eta t} C\int_t^{\infty} e^{-\eta s} W(s)y_0 d s-CW(t)y_0=C\frac{du(t)}{d t}
$$
for any $t \geq 0$,
which implies that
$$
Au(t)=\eta \cdot e^{\eta t} \int_t^{\infty} e^{-\eta s} W(s)y_0 d s-W(t)y_0=\frac{du(t)}{d t},
$$
hence $\frac{du(t)}{d t} \in C\left([0, \infty), S\right)$, i.e., $u \in C^1\left([0, \infty), S\right)$.  Since $u$ and $Au$ are continuous from $[0, \infty)$ to $S$, it follows that
$u \in C\left([0, \infty), [D(A)]\right)$. Thus $u$ is a solution of (6).

Let $t \geq 0$ be fixed,
suppose that $v$ is an arbitrary solution of (6) satisfying the locally $L^0$-Lipschitz condition, define a function $f:[0, t] \rightarrow S$ by $f(s):=C^{2}W(t-s)v(s)$, then it is easy to check that $f$ is $\mathcal{T}_{\varepsilon, \lambda}$-differentiable on $S$ and  $f^{\prime}(s)=0$ for any $s \in [0, t]$.
Since
$$
\begin{aligned}
\| f\left(s_1\right)-f\left(s_2\right)\|= & \left\|C^{2}W\left(t-s_1\right) v\left(s_1\right)-C^{2}W\left(t-s_2\right) v\left(s_2\right)\right\| \\\leq & \left\|C^{2}W\left(t-s_1\right) v\left(s_1\right)-C^{2}W\left(t-s_2\right) v\left(s_1\right)\right\| \\
& +\|C^{2}W\left(t-s_2\right)\left(v\left(s_1\right)-v\left(s_2\right)) \|\right.
\end{aligned}
$$
for any $s_1, s_2\in[0,t]$,
 it follows from Proposition 4.1(5) that $f$ is locally $L^0$-Lipschitz.
Due to Proposition 2.1, $f^{\prime}$ is Riemann integrable and
$$
f(t)-f(0)=\int_0^t f^{\prime}(s) d s=\theta
$$
for any $t \geq 0$, i.e., $C^{3}v(t)=C^{2}W(t)u_0$. Thus $v(t)=W(t)C^{-1}u_0$ for any $t \geq 0$, i.e., the uniqueness of the solution of
(6)  has been proved.
\endproof

Suppose that $\{W(s):s \geq 0\}$ belongs to $ H(M, \xi)$ and $(A, D(A))$ is the infinitesimal generator of $\{W(s):s \geq 0\}$.  For any $y \in C (D(A))$, there exists a $z \in D(A)$ such that $y=C z $. For any $\eta \in G_{\xi}$, let $x=(\eta-A) z $. According to (4), we have
\allowdisplaybreaks
\begin{align*}
R_C(\eta, A) x & =\int_0^{\infty} e^{-\eta s} W(s) x d s \\
& =\int_0^{\infty} e^{-\eta t}W(s)(\eta-A) z d s \\
& =Cz \\
& =y.
\end{align*}
Thus $y \in R(R_C(\eta, A)$ for any $\eta \in G_{\xi}$, i.e., $C(D(A)) \subseteq R(R_C(\eta, A))$.
Due to Theorem 5.1, one can  immediately  obtain Corollary 5.1 below.

\textbf{Corollary 5.1} \textrm{(see \cite{Son})}
Let $\{W(s):s \geq 0\}$ be an exponentially bounded $C$-semigroup on $S$ with the infinitesimal generator $(A, D(A))$. Then, under the
locally $L^0$-Lipschitz condition on the solution, the abstract Cauchy problem
$$
\left\{\begin{array}{l}
\frac{d u(t)}{d t}=A u(t),  \forall t \geq 0,\\
u(0)=u_0 \in C(D(A))
\end{array}\right.
$$
has a unique solution $u(t):=W(t)C^{-1}u_0$ belonging to $C([0, \infty),[D(A)])\cap C^1([0, \infty),S)$.


\begin{thebibliography}{99}


\bibitem{Arendt}
Arendt, W., Batty, C.J.K., Hieber, M. and Neubrander, F., Vector-valued Laplace transforms and Cauchy problems, Birkh\"auser Verlag, Basel, 2001.

\bibitem{Davies}
Davies, E.B. and Pang, M.M.H., The Cauchy problem and a generalization of the Hille-Yosida theorem, \textit{Proc. London Math. Soc.} (3) \textbf{55}, 1987, 181--208.

\bibitem{deLaubenfels1}
deLaubenfels, R., Existence Families, Functional Calculi and Evolution Equations, Springer-Verlag, Berlin, 1994.



\bibitem{D}
Dunford, N. and Schwartz, J.T., Linear Operators (\uppercase\expandafter{\romannumeral 1}), Interscience, New York, 1957.


\bibitem{Guo01}
Guo, T.X., A new approach to random functional analysis, Proceedings of the First China Postdoctoral Academic Conference, The China National Defense and Industry Press, Beijing, 1993, pp. 1150--1154.

\bibitem{TT}
Guo, T.X., Extension theorems of continuous random linear operators on random domains, \textit{J. Math. Anal. Appl.} \textbf{193}(1), 1995, 15--27.

\bibitem{Guo10}
Guo, T.X., Relations between some basic results derived from two kinds of topologies for a random locally convex module, \textit{J. Funct. Anal.} \textbf{258}(9), 2010, 3024--3047.

\bibitem{GuoZhaozeng14}
Guo, T.X., Recent progress in random metric theory and its applications to conditional risk measures, \textit{Sci. China Ser. A} \textbf{54}, 2011, 633--660.

\bibitem{Guo2024}
Guo, T.X., Mu, X.H. and Tu, Q., The relations among the notions of various kinds of stability and their applications, \textit{Banach J. Math. Anal.} \textbf{18}, 2024, 42.

\bibitem{Guo20241}
Guo, T.X., Wang, Y.C., Xu, H.K., Yuan, G. and Chen, G., A noncompact Schauder fixed point theorem in random normed modules and its applications, \textit{Math. Ann.} \textbf{391}, 2025, 3863--3911.

\bibitem{Guo0912}
Guo, T.X., Xiao, H.X. and Chen, X.X., A basic strict separation theorem in random locally convex modules, \textit{Nonlinear Anal.} \textbf{71}(9), 2009, 3794--3804.

\bibitem{GuoZhangwuyang17}
Guo, T.X., Zhang, E.X., Wang, Y.C. and Guo, Z.C., Two fixed point theorems in complete random normed modules and their applications to backward stochastic equations, \textit{J. Math. Anal. Appl.} \textbf{486}(2), 2020, 123644.

\bibitem{Guozhang12}
Guo, T.X. and Zhang, X., Stone's representation theorem of a group of random unitary operators on complete complex random inner product modules, \textit{Sci. Sin. Math.} \textbf{42}(3), 2012, 181--202.

\bibitem{Guozhao}
Guo, T.X. and Zhao, S.E., On the random conjugate spaces of a random locally convex module, \textit{Acta Math. Sin. (Engl. Ser.)} \textbf{28}(4), 2012, 687--696.

\bibitem{Guozhao15}
Guo, T.X., Zhao, S.E. and Zeng, X.L., Random convex analysis (I): Separation and Fenchel-Moreau duality in random locally convex modules, \textit{Sci. Sin. Math.} \textbf{45}(12), 2015, 1961--1980.

\bibitem{Hille48}
Hille, E., Functional Analysis and Semigroups, Amer. Math. Soc. Colloq. Publ., New York, 1948.

\bibitem{LIU2}
Liu, M., Zhang, X. and Dai, L.F., Trotter-Kato approximations of semilinear stochastic evolution equations in Hilbert spaces, \textit{J. Math. Phys.} \textbf{64}(4), 2023, 043506.

\bibitem{LIU1}
Liu, M., Zhang, X. and Dai, L.F., Trotter-Kato approximations of impulsive neutral SPDEs in Hilbert spaces, \textit{Acta Math. Sin. (Engl. Ser.)} \textbf{40}(5), 2024, 1229--1243.


\bibitem{E2019}
Lu\v{c}i\'{c}, D. and Pasqualetto, E., The Serre-Swan theorem for normed modules, \textit{Rend. Circ. Mat. Palermo} \textbf{68}(2), 2019, 385--404.

\bibitem{E}
Lu\v{c}i\'{c}, M., Pasqualetto, E. and Vojnovi\'{c}, I., On the reflexivity properties of Banach bundles and Banach modules, \textit{Banach J. Math. Anal.} \textbf{18}, 2024, 7.

\bibitem{Miyadera86}
Miyadera, I., On the generators of exponentially bounded $C$-semigroups, \textit{Proc. Japan Acad. Ser. A} \textbf{62}, 1986, 239--242.

\bibitem{Phillips1954}
Phillips, R.S., An inversion formula for Laplace transform and semigroups of linear operators, \textit{Ann. Math.} \textbf{59}(2), 1954, 325--356.

\bibitem{Son}
Son, T.C., Thang, D.H. and Oanh, L.T., Exponentially bounded $C$-semigroup and the Cauchy initial value problems in complete random normed modules, \textit{Acta Math. Sin. (Engl. Ser.)} \textbf{40}(9), 2024, 2195--2212.

\bibitem{Thang19}
Thang, D.H., Son, T.C. and Thinh, N., Semigroups of continuous module homomorphisms on complex complete random normed modules, \textit{Lithuanian Math. J.} \textbf{59}(2), 2019, 229--250.

\bibitem{Wu2013}
Wu, M.Z., Fakas' lemma in random locally convex modules and Minkowski-Weyl type results in $L^0(\mathcal{F}, R^n)$, \textit{J. Math. Anal. Appl.} \textbf{404}(2), 2013, 300--309.

\bibitem{WuGuo2015}
Wu, M.Z. and Guo, T.X., A counterexample shows that not every locally $L^0$-convex topology is necessarily induced by a family of $L^0$-seminorms, https://doi.org/10.48550/arXiv.1501.04400, 2015.

\bibitem{WuGuoLong2022}
Wu, M.Z., Guo, T.X. and Long, L., The fundamental theorem of affine geometry in regular $L^0$-modules, \textit{J. Math. Anal. Appl.} \textbf{507}(2), 2022, 125827.

\bibitem{WuZengZhao2022}
Wu, M.Z., Zeng, X.L. and Zhao, S.E., On $L^0$-convex compactness in random locally convex modules, \textit{J. Math. Anal. Appl.} \textbf{515}(2), 2022, 126404.

\bibitem{Yosida1948}
Yosida, K., On the differentiability and the representation of one-parameter semigroups of linear operators, \textit{J. Math. Soc. Japan} \textbf{1}, 1948, 15--21.

\bibitem{Zhangliu13}
Zhang, X. and Liu, M., On almost surely bounded semigroups of random linear operators, \textit{J. Math. Phys.} \textbf{54}(5), 2013, 053517.


\bibitem{Zhangliu25}
Zhang, X. and Liu, M., On the abstract Cauchy problem for evolution equations in a complete random normed module, \textit{J. Math. Phys.} \textbf{66}(1), 2025, 013507.

\bibitem{Zhangliuguo2}
Zhang, X., Liu, M. and Guo, T.X., The Hille-Yosida generation theorem for almost surely bounded $C_0$-semigroups of continuous module homomorphisms, \textit{J. Nonlinear Convex Anal.} \textbf{21}(9), 2020, 1995--2009.

\bibitem{Zhangliuwei}
Zhang, X., Wei, L.L. and Liu, M., $C$-existence families, $C$-semigroups and their associated abstract Cauchy problems in complete random normed modules,
https://doi.org/10.48550/arXiv.2503.03096, 2025.


\end{thebibliography}
\end{document}